\documentclass{amsart}
\newcommand{\R}{{\mathbb R}}
\newcommand{\C}{{\mathbb C}}
\newcommand{\pr}[1]{\left( #1\right)}
\newcommand{\pmd}[1]{\left| #1\right|}
\newcommand{\pq}[1]{\left[ #1\right]}
\newcommand{\sgn}{\,\textnormal{sgn}\,}
\newcommand{\est}[1]{\begin{equation*}\begin{split}#1\end{split}\end{equation*}}
\newcommand{\es}[1]{\begin{equation}\begin{split}#1\end{split}\end{equation}}
\newcommand{\tn}[1]{\textnormal{#1}}
\newtheorem{teo}{Theorem}
\newtheorem{lemma}{Lemma}

\usepackage{amsthm,amsfonts,amssymb,amsmath}
\usepackage{graphicx}

\usepackage{epstopdf}
\usepackage{nicefrac}

\begin{document}

\markboth{Sandro Bettin}
{Second moment of the Zeta-function with unbounded shifts}

\title{The second moment of the Riemann zeta function with unbounded shifts}

\author{SANDRO BETTIN}

\address{Department of Mathematics, University of Bristol, University Walk\\
Bristol, BS8 1TW, United Kingdom}
\email{Sandro.Bettin@bristol.ac.uk}

\maketitle

\begin{abstract}
We prove an asymptotic formula for the second moment (up to height $T$) of the Riemann zeta function with two shifts. The case we deal with is where the real parts of the shifts are very close to zero and the imaginary parts can grow up to $T^{2-\varepsilon}$, for any $\varepsilon>0.$ 
\end{abstract}

\section{Introduction}

An important problem in analytic number theory is to understand the moments of the Riemann zeta function
\est{
I_k(T)=\int_0^T\pmd{\zeta\pr{\frac12+it}}^{2k}\,\tn{d}t.
}
The knowledge of the asymptotic behavior of $I_k(T)$ would give important information about the maximal order of the Riemann zeta function on the critical line and about the zeros of this function and therefore about the distribution of prime numbers. Unfortunately, the asymptotic is known just for $k=1$ and $k=2$. Specifically, in 1918 Hardy and Littlewood~\cite{HL} proved
\est{
I_1(T)\sim T\log T
}
and in 1926 Ingham~\cite{ING} proved
\est{
I_2(T)\sim \frac1{2\pi^2}T\log^4 T.
}
For other values of $k$ the problem is still open. For positive real numbers $k$ it is conjectured that 
\est{
I_k(T)\sim C_kT\log^{k^2}T
}
for a positive constant $C_k$.
Conrey and Ghosh~\cite{CGH} conjectured a value for $C_6$ and later Conrey and Gonek~\cite{CGK} made a conjecture for $C_8$. In 2000, using results on random matrix theory, Keating and Snaith~\cite{KS} conjectured the value of $C_k$ for all positive $k$.

It is natural to consider also the shifted moments
\est{
\int_0^T\prod_{a\in A}\zeta\pr{\frac12+a+it}\prod_{b\in B}\zeta\pr{\frac12-b-it}\,\tn{d}t
}
with $A$, $B$ finite sets of complex numbers with the same cardinality. In~\cite{ING}, Ingham dealt with the case with just two zeta functions, obtaining that the second shifted moment
\est{
\int_\kappa^T\zeta\pr{\frac12+a+it}\zeta\pr{\frac12-b-it}\,\tn{d}t,
}
is asymptotic to
\est{
\int_{1}^T\pr{\zeta(1+c)+\tau^{-c}\zeta(1-c)}\,\tn{d}\tau,
}
where $c=a-b$ and provided that $a$ and $b$ are bounded and $\kappa>1+\max\pr{|a|,|b|}$. 

For some purposes it is useful to have results also for unbounded shifts. For example, in~\cite{CS} the computation of the one-level density for zeros of quadratic Dirichlet $L$-functions is obtained using the ratios conjecture and integrating over the shifts on an unbounded interval. The ratios conjecture suggests a formula for averages of ratios of products of zeta functions (of which the second moment above is a special case). This conjecture has been very useful, but it is completely unknown on which range of shift parameters it holds.

As in Ingham's work, in this paper we deal with the second shifted moment, but allowing the imaginary part of the shifts $a$ and $b$ to grow with $T.$ Essentially, we obtain that Ingham's asymptotic formula is still true if the imaginary parts of $a$ and $b$ are less than $T^{2-\varepsilon}$. It is very likely that the formula still holds up to height $T^{4-\varepsilon}$, though it might not be possible to exceed this limit because of omega results in the error term for $I_1(T)$.

This work lends support to the ratios conjecture formula holding over a large range of shift parameters.

For the sake of convenience, we assume that the real part of $a$ and $b$ are close to 0, but the same method should work also for small bounded values of $\Re(a)$ and $\Re(b)$. More specifically, we obtain the following.
\begin{teo}
Let $T\geq2$, $a=a(T),$ $b=b(T)\in\C$ such that
\est{
\Re(a)&\ll\frac1{\log T},\\
\Re(b)&\ll\frac1{\log T}\\
}
and let $c=a-b.$ Then
\est{
\int_0^T\zeta\pr{\frac12+a+it}\zeta&\pr{\frac12-b-it}\,\tn{d}t\\
=&\int_{0}^{T}\pr{\zeta(1+c)+\zeta(1-c)\chi\pr{\frac12+a+it}\chi\pr{\frac12-b-it}}\,\tn{d}t+\\
&+O\pr{\pr{T^{\frac 12}+\sqrt{|\Im(a)|}+\sqrt{|\Im(b)|}}\log^2 T},\\
}
as $T\rightarrow\infty.$ If $c=0,$ then the integrand on the right hand side has to be interpreted as the pointwise limit for $c\rightarrow0$.
\end{teo}
The $\chi$ function in the theorem is the function that appears in the functional equation, 
\est{
\zeta(1-s)=\chi(1-s)\zeta(s),
}
that is
\est{
\chi(1-s)=2(2\pi)^{-s}\Gamma(s)\cos\frac{\pi s}2.
}

The ideas of the proof stem from the proof of Theorem 7.4 in~\cite{TITCH}, about the second moment of the Riemann Zeta function without shifts. Before giving the proof, we prove two lemmas, the first is the analogue with shift of Theorem 12.2 in~\cite{TITCH}, the second gives an asymptotic formula for $\chi\pr{\frac12+a+it}\chi\pr{\frac12-b-it}$.

\section{Proof of the Theorem}

In the whole article we consider $a$, $b$ and $c=a-b$ to be complex numbers and we write
\est{
a&=\alpha+i\alpha',\\
b&=\beta+i\beta',\\
c&=\gamma+i\gamma',
}
with $\alpha,\alpha',\beta,\beta',\gamma,\gamma'\in\R.$
\begin{lemma}\label{a}
Let $x\geq2$ and $c=c(x)$. Let
\est{
D_c(x)=\sum_{mn\leq x}\frac 1{n^c}=\sum_{n\leq x}\sigma_{-c}(n),
}
where
\est{
\sigma_q(n)=\sum_{d|n} d^q.
}
Assume
\est{|\gamma|\ll\frac1{\log x}.} 
Then, for any $\varepsilon>0$, we have
\est{
D_c(x)&=\zeta(1+c)x+\zeta(1-c)\frac{x^{1-c}}{1-c}+O\pr{x^{\frac13+\varepsilon}+|\gamma'|^{\frac12}\log^2 x}\\
&=\int_1^x\pr{\zeta(1+c)+\zeta(1-c)u^{-c}}\,\tn{d}u+O\pr{x^{\frac13+\varepsilon}+|\gamma'|^{\frac12}\log^2x},
}
as $x\rightarrow\infty.$ As above, if $c=0$ the right hand side has to be interpreted as the limit for $c\rightarrow0.$
\end{lemma}
\begin{proof}
Clearly, we can assume $x$ is half an odd integer. Let $Q\geq1,$ $\varepsilon>0$ and $\eta=\frac1{\log x}+\max(1,1-\gamma)$. Applying Lemma 3.12 in~\cite{TITCH} with $a_n=\sigma_{-c}(n)$, $\psi(n)=n^\varepsilon$ and with the $c$ of the lemma equal to $\eta$, we find that
\es{\label{pr}
D_c(x)&=\frac1{2\pi i}\int_{\eta-iQ}^{\eta+iQ}\!\!\zeta(s)\zeta(s+c)\frac{x^s}{s}\,\tn{d}s\!\!+O\pr{\frac{x^{\eta}}{Q}\log^2x}\!+O\pr{\frac{x^{1+\varepsilon+\max\pr{-\gamma,0}}\log x}{Q}}\\
&=\frac1{2\pi i}\int_{\eta-iQ}^{\eta+iQ}\zeta(s)\zeta(s+c)\frac{x^s}{s}\,\tn{d}s+O\pr{\frac{x^{1+\varepsilon}\log x}{Q}},
}
where we used (1.3.1) in~\cite{TITCH}. Take $Q=x^\frac23$ (replacing $Q$ with $Q+1$ if $\pmd{x^\frac23-|\gamma'|}<\frac12$) and move the path of integration to $\nu=-\frac1{\log x}+\min\pr{0,-\gamma}.$ Since the integrand is $O\pr{x^{\frac13}\log^2x+|\gamma'|^{\frac12}x^{-\frac13}}$ on the horizontal lines, we find 
\es{\label{a1}
\frac1{2\pi i}\int_{\eta-iQ}^{\eta+iQ}\zeta(s)\zeta(s+c)\frac{x^s}{s}\,\tn{d}s=&\frac1{2\pi i}\int_{\nu-iQ}^{\nu+iQ}\zeta(s)\zeta(s+c)\frac{x^s}{s}\,\tn{d}s+\zeta(1+c)x+\\
&\!+\zeta(1-c)\frac{x^{1-c}}{1-c}\!+\zeta(c)\!+O\pr{x^{\frac13}\log^2x+|\gamma'|^{\frac12}x^{-\frac13}},
}
where $\zeta(c)$ has to be omitted if $|\gamma'|\geq x^\frac23+\frac12$ and so it can be always inserted in the error term, since 
\est{
\zeta(c)\ll|c|^{\frac12}\log\pr{|c|+2},
}
by the functional equation and Theorem 3.5 in~\cite{TITCH}. Moreover, expanding the zetas as a series, we find 
\es{\label{a2}
\int_{\nu-iQ}^{\nu+iQ}\zeta(s)\zeta(s+c)\frac{x^s}{s}\,\tn{d}s&=\int_{\nu-iQ}^{\nu+iQ}\chi(s)\chi(s+c)\zeta(1-s)\zeta(1-s-c)\frac{x^s}{s}\,\tn{d}s\\
&=ix^{\nu}\sum_{n=1}^{\infty}\frac{\sigma_{-c} (n)}{n^{1-c-\nu}}\int_{-Q}^{Q}\chi(\nu+it)\chi(\nu+c+it)\frac{(xn)^{it}}{\nu+it}\,\tn{d}t.
}
The asymptotic expansion in a strip for $\chi(s)$ (see, for example, (4.12.3) in~\cite{TITCH}) gives
\est{
i\chi(\nu+it)\chi(\nu+&c+it)\frac{(xn)^{it}}{\nu+it}=\\
=&\pr{\frac{|t|}{2\pi}}^{\frac12-\nu-it}\pr{\frac{|t+\gamma'|}{2\pi}}^{\frac12-\nu-c-it}\times\\
&\times e^{i\pr{2t+\gamma'+\frac{\pi}{4}(\sgn (t)+\sgn(t+\gamma'))}}\frac{(xn)^{it}}{t}\pr{1+O\pr{\frac1{|t|}+\frac1{|t+\gamma'|}}}\\
=&\,g(t)e^{ih(t)}\pr{1+O\pr{\frac1{|t|}+\frac1{|t+\gamma'|}}},
}
where
\est{
g(t)=&\frac1t\pr{\frac{|t|}{2\pi}}^{\frac12-\nu}\pr{\frac{|t+\gamma'|}{2\pi}}^{\frac12-\nu-\gamma},\\
h(t)=&-t\log\pr{\frac{|t|}{2\pi}}-\pr{\gamma'+t}\log\pr{\frac{|t+\gamma'|}{2\pi}}+2t+\gamma'+\\
&+\frac{\pi}{4}(\sgn(t)+\sgn(t+\gamma'))+t\log (xn),\\
h'(t)=&-\log\pr{\frac{|t|}{2\pi}}-\log\pr{\frac{|t+\gamma'|}{2\pi}}+\log (xn)=-\log\pr{\frac{|t+\gamma'||t|}{(2\pi)^2xn}},\\
h''(t)=&-\frac{2t+\gamma'}{t(t+\gamma')}.
}
Firstly let's consider the case $|\gamma'|< 3Q$. If $n>\frac{Q^2}x$, we have trivially that $h'\gg1$. If $n\leq\frac{Q^2}x$, 
defining
\est{
V&=\left\{t\mid \left||t+\gamma'||t|-(2\pi)^2xn\right|<\sqrt nx^\frac56\right\},\\
U&=[-Q,Q]\setminus V,
} 
we have
\est{
h'(t)\gg\frac{1}{x^\frac16n^\frac12},
}
for $t\in U$,
and
\est{
\tn m\pr{V}\ll x^\frac13,
}
where $\tn m\pr{V}$ is the measure of $V$.
Thus, using Lemma 4.3 in~\cite{TITCH},  we have that
\es{\label{a3}
ix^\nu\sum_{n=1}^{\infty}\frac{\sigma_{-c} (n)}{n^{1-c-\nu}}&\int_{-Q}^{Q}\chi(\nu+it)\chi(\nu+c+it)\frac{(xn)^{it}}{\nu+it}\,\tn{d}t=\\
&=x^\nu\sum_{n\leq\frac{Q^2}x}\frac{\sigma_{-c} (n)}{n^{1-c-\nu}}\int_{-Q}^{Q}g(t)e^{ih(t)}\,\tn{d}t+O\pr{\log^2 x}\\
&=x^\nu\sum_{n\leq\frac{Q^2}x}\frac{\sigma_{-c} (n)}{n^{1-c-\nu}}\int_Ug(t)e^{ih(t)}\,\tn{d}t+O\pr{x^\frac13\log^2 x}\\
&\ll x^\frac13\log^2 x,
}
if $|\gamma'|< 3Q$.
Now, if $|\gamma'|\geq 3Q,$ we have
\est{
h''(t)\gg\frac1Q,
}
in $[-Q,Q]$. Thus, using Lemma 4.5 in~\cite{TITCH}, we have
\est{
ix^\nu\sum_{n=1}^{\infty}\frac{\sigma_{-c} (n)}{n^{1-c-\nu}}&\int_{-Q}^{Q}\chi(\nu+it)\chi(\nu+c+it)\frac{(xn)^{it}}{\nu+it}\,\tn{d}t\ll |\gamma'|^\frac12\log^2 x.
}
This equation, together with~\eqref{pr}-\eqref{a3}, implies the stated result.
\end{proof}
From now on, let $T\geq2$ and assume $a$ and $b$ are functions of $T$ such that
\est{
|\alpha|,|\beta|&\ll\frac1{\log T}\\
|\alpha'|,|\beta'|&\ll T^2.
}
Moreover, let $u$ and $v$ such that
\es{\label{uv}
T^u&=\max(|\alpha'|,T),\\
T^v&=\max(|\beta'|,T).
}
\begin{lemma}\label{chi}
Let  
\est{
|t+\alpha'|&>10|c|,\\
|t+\beta'|&>10|c|.
} 
Then
\est{
\chi\pr{\frac12+a+it}&\chi\pr{\frac12-b-it}=\\
&e^{-c\log \frac{|t+\beta'|}{2\pi}+\pr{a+it}\log\pr{1+\frac{c}{\frac12-a-it}}+c}\pr{1+O\pr{\frac{|c|}{|t+\alpha'|}}}.
}
\end{lemma}
\begin{proof}
Since $\chi(s)\chi(1-s)=1,$ we have
\es{\label{b1}
\chi\pr{\frac12+a+it}\chi\pr{\frac12-b-it}&=\frac{\chi\pr{\frac12+a+it}}{\chi\pr{\frac12+b+it}}\\
&=(2\pi)^{c}\frac{\cos\pr{\frac\pi2\pr{\frac12-a-it}}}{\cos\pr{\frac\pi2\pr{\frac12-b-it}}}\frac{\Gamma\pr{\frac12-a-it}}{\Gamma\pr{\frac12-b-it}}.
}
Stirling's formula, as expressed in~\cite{RAD}, (21.1), states
\est{
\log\Gamma(s)=\left(s-\frac12\right)\log s-s+\frac12\log 2\pi+R(s),
}
where
\est{
R(s)=\int_0^{+\infty}\frac{g(x)}{(x+s)^2}\,\tn{d}x
}
and $g(x)=\{x\}(\{x\}-1)/2$. Therefore
\est{
\frac{\Gamma\pr{\frac12-a-it}}{\Gamma\pr{\frac12-b-it}}=e^{-\pr{a+it}\log\pr{\frac12-a-it}+\pr{b+it}\log\pr{\frac12-b-it}+c+\Delta_{a,b}\pr{t}},
}
where
\est{\Delta_{a,b}\pr{t}=R\pr{\frac12-a-it}-R\pr{\frac12-b-it}.}
Thus
\es{\label{b2}
\frac{\Gamma\pr{\frac12-a-it}}{\Gamma\pr{\frac12-b-it}}=&\exp\Bigg(\pr{a+it}\log\pr{1+\frac{c}{\frac12-a-it}}-c\log\pr{\frac12-b-it}+\\
&+c+\Delta_{a,b}(t)\Bigg)\\
=&\exp\Bigg(\pr{a+it}\log\pr{1+\frac{c}{\frac12-a-it}}+c-c\log|t+\beta'|+\\
&+\Delta_{a,b}(t)-\frac{\pi i c}2\sgn\pr{t+\beta'}+O\pr{\frac{|c|}{|\alpha'+t|}}\Bigg)\\
}
Moreover,
\es{\label{b3}
\frac{\cos\pr{\frac\pi2\pr{\frac12-a-it}}}{\cos\pr{\frac\pi2\pr{\frac12-b-it}}}&=\frac{\cos\pr{\frac\pi2\pr{\frac12-b-it}-\frac{\pi c}2}}{\cos\pr{\frac\pi2\pr{\frac12-b-it}}}\\
&=\cos{\frac{\pi c}2}+\sin{\frac{\pi c}2}\tan\pr{\frac\pi2\pr{\frac12-b-it}}\\
&=\cos{\frac{\pi c}2}+\sin{\frac{\pi c}2}\pr{i\sgn\pr{t+\beta'}+O\pr{e^{-\pi|t+\beta'|}}}\\
&=e^{\frac{\pi i c}2\sgn\pr{t+\beta'}}+O\pr{\pmd{\sin\pr{\frac {\pi c}2}}e^{-\pi|t+\beta'|}}\\
&=e^{\frac{\pi i c}2\sgn\pr{t+\beta'}}+O\pr{|c|e^{-\frac\pi2|t+\beta'|}}\\
&=e^{\frac{\pi i c}2\sgn\pr{t+\beta'}}+O\pr{\frac{|c|}{|t+\alpha'|}}.
}
Finally, 
\est{
|\Delta_{a,b}(t)|&=\left|\int_0^{+\infty}\frac{g(x)}{(x+\frac12-a-it)^2}-\frac{g(x)}{(x+\frac12-b-it)^2}\,\tn{d}x\right|\\
&\ll \int_0^{+\infty}\left|\frac{1}{(x+\frac12-a-it)^2}-\frac{1}{(x+\frac12-b-it)^2}\right|\,\tn{d}x\\
&\ll \int_0^{+\infty}\left|\frac{c^2+2c\pr{x+\frac12-a-it}}{(x+\frac12-a-it)^2(x+\frac12-b-it)^2}\right|\,\tn{d}x\\
&\ll \frac{|c|^2+2|c||t+\alpha'|}{|t+\alpha'|^2}\int_0^{+\infty}\frac{1}{\pmd{x+\frac12-b-it}^2}\,\tn{d}x\\
&\ll \frac{|c|^2+2|c||t+\alpha'|}{|t+\alpha'|^2\pr{|t+\alpha'|+1|}}.
}
This equation, together with~\eqref{b1}-\eqref{b3}, gives the result.
\end{proof}
\begin{proof}[Proof of Theorem 1]
We have
\est{
\int_0^T\zeta\pr{\frac12+a+it}\zeta\pr{\frac12-b-it}\,\tn{d}t&=-i\int_{\frac12}^{\frac12+iT}\zeta(s+a)\zeta(1-s-b)\,\tn{d}s\\
&=-i\int_{\frac12}^{\frac12+iT}\zeta(s+a)\chi(1-s-b)\zeta(s+b)\,\tn{d}s.
}
Moving the path of integration to $\delta=\varepsilon+\max\pr{1-\alpha,1-\beta}$, with $\varepsilon=\frac1{\log T}$, we obtain
\est{
\int_0^T\zeta\pr{\frac12+a+it}\zeta\pr{\frac12-b-it}\,\tn{d}t=&-i\int_{\delta}^{\delta+iT}\!\!\!\zeta(s+a)\chi(1-s-b)\zeta(s+b)\,\tn{d}s+\\
&+O\pr{T^{r(u+v)}+T^{\frac v2}\log^2T},
}
where $u$ and $v$ are as defined in~\eqref{uv} and $r$ is any constant such that $\zeta\pr{t}\ll t^{r-\epsilon}$, for some $\epsilon>0$.
Now we have
\est{
\int_{\delta}^{\delta+iT}\zeta(s+a)\chi(1-s-b)\zeta(s+b)\,\tn{d}s=&\int_{\delta+b}^{\delta+b+iT}\zeta(s+c)\chi(1-s)\zeta(s)\,\tn{d}s\\
=&\sum_{n,m}\frac{1}{n^c}\int_{\delta+b}^{\delta+b+iT}\chi(1-s)(nm)^{-s}\,\tn{d}s\\
}
and, if $\beta'>0$, by (7.4.2) and (7.4.3) in~\cite{TITCH}, this is
\est{
\sum_{n,m}\frac{1}{n^c}\int_{\delta+b}^{\delta+b+iT}&\chi(1-s)(nm)^{-s}\,\tn{d}s\\
=&\,2\pi i\sum_{\frac{\beta'}{2\pi}<nm\leq\frac{T+\beta'}{2\pi}}\frac{1}{n^c}+O\pr{T^{\frac v2}}+\\
&+O\pr{\sum_{\substack{|2\pi mn-T|>\frac12,\\|2\pi mn-T-\beta'|>\frac12}}\frac{1}{mn}\frac{T^{\frac v2}}{\min\pr{1,\pmd{\log \frac{|\beta'|}{2\pi mn}},\pmd{\log \frac{|T+\beta'|}{2\pi mn}}}}}\\
=&\,2\pi i\sum_{\frac{\beta'}{2\pi}<nm\leq\frac{T+\beta'}{2\pi}}\frac{1}{n^c}+O\pr{T^{\frac v2}\log^2 T},
}
since
\est{
\sum_{|mn-Q|>\frac12}\frac{1}{\pmd{\log \frac{Q}{mn}}mn}\ll&\sum_{m,n}\frac{1}{mn}+\sum_{\substack{|mn-Q|>\frac12,\\ \frac Q{2m}<n<\frac{2Q}m}}\frac{1}{\pmd{\log \frac{Q}{ mn}}mn}\\
\ll&\log^2Q+\sum_{m,r\ll Q}\frac{1}{mr}\ll\log^2 Q.
}
If $\beta'\leq0$ we can proceed in the same way, but we have to replace 
\est{
\sum_{\frac{\beta'}{2\pi}<nm\leq\frac{T+\beta'}{2\pi}}\frac{1}{n^c}
}
with
\begin{equation}
\begin{cases}
\sum_{0<nm\leq\frac{T+\beta'}{2\pi}}\frac{1}{n^c} + \sum_{0<nm\leq\frac{-\beta'}{2\pi}}\frac{1}{n^c} & \text{if $-T\leq\beta'\leq 0 $,}\\
\sum_{-\frac{T+\beta'}{2\pi}\leq nm<\frac{-\beta'}{2\pi}}\frac{1}{n^c} & \text{if $\beta'\leq-T $.}
\end{cases}
\end{equation}
Therefore, by Lemma~\ref{a}, we have
\est{
\int_0^T\zeta\pr{\frac12+a+it}\zeta\pr{\frac12-b-it}&\,\tn{d}t=2\pi\int_{\frac{\beta'}{2\pi}}^{\frac{T+\beta'}{2\pi}}\pr{\zeta(1+c)+\zeta(1-c)|x|^{-c}}\,\tn{d}x+\\
&+O\pr{|\gamma'|^{\frac12}\log^2 T+T^{r\pr{u+v}}+T^{\frac v2}\log^2T}\\
=&\int_{0}^{T}\pr{\zeta(1+c)+\zeta(1-c)\pr{\frac{|t+\beta'|}{2\pi}}^{-c}}\,\tn{d}t+\\
&+O\pr{T^{\frac v2}\log^2T+T^\frac u2\log^2 T},
}
taking, for example, $r=\frac14.$\par
If $\gamma'\gg 1$, Lemma 4.3 in~\cite{TITCH} gives 
\est{
\int_0^T\chi\pr{\frac12+a+it}\chi\pr{\frac12-b-it}\,\tn{d}t\ll \frac{T^{v}}{|\gamma'|}
}
and clearly
\est{
\int_0^T\pr{\frac{|t+\beta'|}{2\pi}}^{-c}\,\tn{d}t\ll \frac{T^{v}}{|\gamma'|}.
}
Therefore, if $|\gamma'|>\frac{T^\frac v2}{\log^2 T},$  we have
\es{\label{ff}
\int_0^T\zeta\pr{\frac12+a+it}&\zeta\pr{\frac12-b-it}\,\tn{d}t=\\
=&\int_{0}^{T}\pr{\zeta(1+c)+\zeta(1-c)\chi\pr{\frac12+a+it}\chi\pr{\frac12-b-it}}\,\tn{d}t+\\
&+O\pr{T^{\frac v2}\log^2 T+T^{\frac u 2}\log^2T}\\
}
Finally, assume $|\gamma'|\leq\frac{T^\frac v2}{\log^2 T}$ (and so also $\gamma'\ll\frac{T^\frac u2}{\log^2 T}$) and let
\est{
W_1&=[0,T]\cap\big\{t\mid|t+\alpha'|,|t+\beta'|>10|c|\big\},\\
W_2&=[0,T]\cap\big\{t\mid|t+\alpha'|,|t+\beta'|\leq10|c|\big\}.
} 
We have
\est{
\int_{W_2}\chi\pr{\frac12+a+it}\chi\pr{\frac12-b-it}\,\tn{d}t&\ll|c|\ll\frac{T^\frac v2}{\log^2 T}
}
and so
\est{
\int_{W_2}\chi\pr{\frac12+a+it}\chi\pr{\frac12-b-it}\,\tn{d}t=\int_{W_2}\pr{\frac{|t+\beta'|}{2\pi}}^{-c}\,\tn{d}t+O\pr{\frac{T^\frac v2}{\log^2 T}}.
}
Furthermore, by Lemma~\ref{chi}, we have
\est{
\int_{W_1}\chi\pr{\frac12+a+it}\chi\pr{\frac12-b-it}\,\tn{d}t=&\int_{W_1}e^{-c\log \frac{|t+\beta'|}{2\pi}+\pr{a+it}\log\pr{1+\frac{c}{\frac12-a-it}}+c}\,\tn{d}t\\
&+O\pr{|c|\log T}
}
and
\est{
\int_{W_1}e^{-c\log \frac{|t+\beta'|}{2\pi}}&\pr{e^{\pr{a+it}\log\pr{1+\frac{c}{\frac12-a-it}}+c}-1}\,\tn{d}t=\\
=&\int_{W_1}\frac{1}{1-c}e^{(1-c)\log \frac{|t+\beta'|}{2\pi}}\Bigg(i\log\pr{1+\frac{c}{\frac12-a-it}}+\\
&+\frac{ic(a+it)}{\pr{\frac12-a-it}\pr{\frac12-b-it}}\Bigg)e^{\pr{a+it}\log\pr{1+\frac{c}{\frac12-a-it}}+c}\,\tn{d}t+\\
&+\pq{\frac{1}{1-c}e^{(1-c)\log \frac{|t+\beta'|}{2\pi}}\pr{e^{\pr{a+it}\log\pr{1+\frac{c}{\frac12-a-it}}+c}-1}}_{\partial W_1}\\
=&\int_{W_1}\frac{1}{1-c}e^{(1-c)\log \frac{|t+\beta'|}{2\pi}}\pr{O\pr{\frac{|c|^2}{|t+\alpha'|^2}}}\,\tn{d}t+\\
&+\pq{\frac{1}{1-c}e^{(1-c)\log \frac{|t+\beta'|}{2\pi}}\pr{O\pr{\frac{|c^2|}{|t+\alpha'|}}}}_{\partial W_1}\\
\ll&|c|\log T\ll\frac{T^{\frac v2}}{\log T}
}
Thus~\eqref{ff} holds also for  $|\gamma'|\leq\frac{T^\frac v2}{\log^2 T}$ and so the proof is complete.
\end{proof}

\section*{Acknowledgments}

I am very grateful to Brian Conrey and Nina Snaith for their help and encouragement.

Research was partially supported by the ``Fondazione Ing. Aldo Gini''. This paper was started when the author was visiting the American Institute of Mathematics (AIM).
\appendix

\end{document}